\documentclass[11pt,letterpaper]{article}
\usepackage[margin=1in]{geometry}
\usepackage[utf8]{inputenc}
\usepackage{amssymb,amsmath,amsthm}
\usepackage{mathrsfs}
\usepackage{hyperref}
\usepackage{cleveref}
\usepackage{graphicx}

\newtheorem{theorem}{Theorem}[section]
\newtheorem{lemma}[theorem]{Lemma}
\newtheorem{proposition}[theorem]{Proposition}
\newtheorem{corollary}[theorem]{Corollary}

\theoremstyle{definition}
\newtheorem{definition}[theorem]{Definition}
\newtheorem{example}[theorem]{Example}
\newtheorem*{pf}{Proof}

\theoremstyle{remark}
\newtheorem*{remark}{Remark}

\crefname{theorem}{Theorem}{Theorems}
\crefname{lemma}{Lemma}{Lemmas}
\crefname{proposition}{Proposition}{Propositions}
\crefname{corollary}{Corollary}{Corollaries}
\crefname{definition}{Definition}{Definitions}
\crefname{example}{Example}{Examples}
\crefname{remark}{Remark}{Remarks}
\crefname{figure}{Figure}{Figures}

\newcommand{\compl}[1]{{#1}^{\mathsf{c}}} %set complement

\newcommand{\cntrlmapf}[3]{G_{{#1},{#2};{#3}}}
\newcommand{\cntrlmapp}{G}
\newcommand{\domattr}[1]{\mathcal{D}({#1})}
\newcommand{\fnorm}[1]{\| #1\|_\text{sup}}
\newcommand{\flow}[2]{\varphi({#1},{#2})}
\newcommand{\gflow}[3]{\varphi({#1},{#2},{#3})}
\newcommand{\gflowf}[4]{\varphi({#1},{#2},{#3};{#4})}
\newcommand{\gflowc}[3]{\varphi^{#1}_{#3}({#2})}
\newcommand{\gspaceI}[1]{L^\infty({#1},\R^n)}
\newcommand{\Infnorm}[1]{\left\| #1 \right\|_\infty} %supnorm
\newcommand{\intensity}[1]{\mu({#1})}
\newcommand{\norm}[1]{\left\| #1 \right\|} % norm
\newcommand{\R}{\mathbb{R}}                % real numbers
\newcommand{\reach}[2]{P_{#1}({#2})}  % reachable set
\newcommand{\Reach}[2]{\overline{P_{#1}({#2})}} %closure of reachable set

\begin{document}

\title{Intensity---A Metric Approach to Quantifying \\ Attractor Robustness in ODEs}
\author{Katherine J. Meyer\footnote{Department of Mathematics and Statistics, Carleton College, Northfield MN. \textit{kjmeyer@carleton.edu}} and Richard P. McGehee\footnote{School of Mathematics, University of Minnesota, Minnepolis MN. \textit{mcgehee@umn.edu}}}
\date{}
\maketitle

\vspace{-0.8cm}

\begin{abstract}
Although mathematical models do not fully match reality, robustness of dynamical objects to perturbation helps bridge from theoretical to real-world dynamical systems. Classical theories of structural stability and isolated invariant sets treat robustness of qualitative dynamics to sufficiently small errors. But they do not indicate just how large a perturbation can become before the qualitative behavior of our system changes fundamentally. Here we introduce a quantity, intensity of attraction, that measures the robustness of attractors in metric terms. Working in the setting of ordinary differential equations on $\mathbb{R}^n$, we consider robustness to vector field perturbations that are time-dependent or -independent. We define intensity in a control-theoretic framework, based on the magnitude of control needed to steer trajectories out of a domain of attraction. Our main result is that intensity also quantifies the robustness of an attractor to time-independent vector field perturbations; we prove this by connecting the reachable sets of control theory to isolating blocks of Conley theory. 
In addition to treating classical questions of robustness in a new metric framework, intensity of attraction offers a novel tool for resilience quantification in ecological applications. Unlike many measurements of resilience, intensity detects the strength of transient dynamics in a domain of attraction.
\end{abstract}


\section{Introduction}

Although differential equations are imperfect models of reality, they capture the essential dynamics of many real-world systems remarkably well. This utility of dynamical systems theory stems in part from studying robust dynamical objects. For example, classical results about structurally stable systems \cite{andronov1937rough,peixoto1962structural,smale1967differentiable} and isolated invariant sets \cite{conleyeaston1971, conley1978} guarantee that qualitative dynamics persist through sufficiently small errors, under suitable conditions of nondegeneracy. In this work we push beyond ``sufficiently small"  to explore just how wrong our models can be while still providing meaningful information about the long-term behavior of a system. More precisely, we ask in metric terms how different a vector field can become  while retaining---in some sense---an original attractor.

To answer this question, we adapt the concept of intensity of attraction from the setting of maps on compact metric spaces \cite{mcgehee1988} to the setting of vector fields on $\R^n$. We define the intensity of an attractor in \cref{sec:IoA}. What makes intensity so well-suited to our purpose is that it  
\begin{itemize}
\item[(i)] quantifies how far an attractor persists in the context of both time-independent and time-varying changes to a vector field $f:\R^n\to\R^n$,
\item[(ii)] can be computed, and 
\item[(iii)] is well-defined for ODE systems of any finite dimension.
\end{itemize}
Our main result, \cref{main}, links the time-dependent and time-independent vector field perturbations in (i): roughly speaking, if we define the intensity $\mu$ of an attractor $A$ for a vector field $f$ as the magnitude of time-varying control needed to escape $A$'s domain of attraction, then the attractor continues in a related form through any time-\emph{independent} perturbations to $f$ that stay below $\mu$ in an appropriate metric.
This connection reduces the problem of attractor persistence in an infinite-dimensional space to the well-studied problem of computing reachable sets in a control framework \cite{BeynRieger2007,Rieger2015, Rieger2016, Sandberg2008,Rasmussen2017,Baier2013,deWeerdt2011}.

In addition to quantifying tolerable modeling errors, intensity of attraction provides a new tool for resilience quantification for applications in ecology and other fields \cite{carpenter2001,meyer2016,meyer2018}. In particular, intensity measures the resilience of a particular regime---modeled as a domain of attraction---to continuous-time, exogeneous disturbances.  Unlike common resilience metrics like eigenvalues or basin size \cite{meyer2016} that are based solely on invariant sets, intensity detects the strength of attracting dynamics over transient regions of state space.

We illustrate some of the preceeding points about intensity of attraction in a simple example before outlining the paper's full contents.

\begin{example}\label{ex:intro}
Suppose the differential equations $x'=f(x)$ and $y'=\widehat f(y)$ model the dynamics of biological populations $x$ and $y$. The vector fields $f:\R\to\R$ and $\widehat f:\R\to\R$ graphed in \cref{introfig}(a)  generate topologically equivalent dynamical systems, each with a repelling equilibrium at the origin and an attracting equilibrium at $1$. In the absence of error or disturbance, we expect each biological population to approach a carrying capacity at $1$. But this topological picture misses a metric difference between the two models: $f$ lies above $\widehat f$ over the transient interval $(0,1)$ of state space. The intensity of each attractor $\{x=1\}$ and $\{y=1\}$ captures this metric difference: in this case, intensity boils down to the maximum values $\mu$ and $\widehat \mu$ that the functions $f$ and $\widehat f$ acheive, respectively, over the interval $(0,1)$. 

\begin{figure}[h]
\centering
\includegraphics[width=\textwidth]{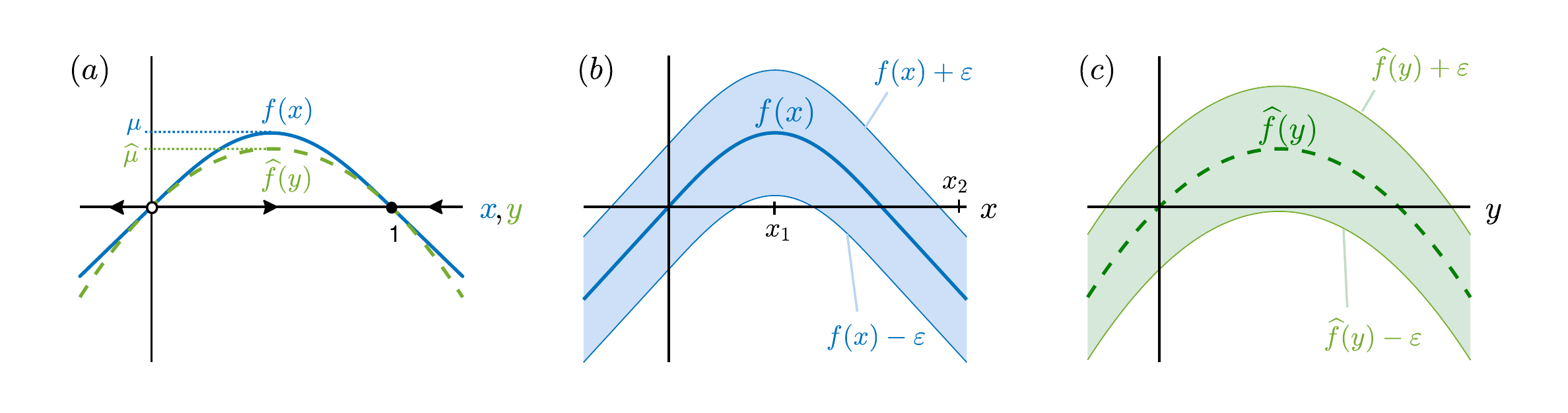}
\vspace*{-1cm} 
\caption{Motivating example. (a) Vector fields $f$ and $\widehat f$ generate topologically equivalent flows on $\R$, illustrated by the horizontal phase line. Intensities $\mu$ and $\widehat\mu$ of the mutual attractor at $1$ capture metric differences between the dynamics. (b) If $f(x)$ is known to an error tolerance of $\epsilon$ with $0<\epsilon<\mu$ then we have $f(x_1)\pm\epsilon>0$ and $f(x_2)\pm\epsilon<0$, suggesting that some trajectories $x(t)$ remain bounded and positive. (c) Under the same error tolerance $\epsilon$ on $\widehat f(y)$, we have $\widehat f(y)-\epsilon <0$ for all $y\in \R$, so it is possible that no trajectories stay positive. }
\label{introfig}
\end{figure}

The fact that the intensity $\mu$ of $\{x=1\}$ exceeds the intensity $\widehat\mu$ of $\{y=1\}$ holds implications for both population resilience and model interpretation.
If we take the models at face value, then $\mu>\widehat\mu$ suggests that population $x$ is more resilient to external time-dependent forcings (such as climate variability or seasonal harvests, as in \cite{meyer2018}) than  population $y$, based on its higher rebound rate at depressed population levels. 
 Notably, two common resilience quantifiers based on invariant sets fail to detect this difference: $f$ and $\widehat f$ share the same distance between equilibria and the same eigenvalues (slopes) at equilibria. 
 
On the other hand, suppose models $x'= f(x)$ and $y'=\widehat f(y)$ originated from data on real-world populations $x$ and $y$, and the values of $f$ and $\widehat f$ are known to an accuracy of $\varepsilon$, where $\widehat\mu<\varepsilon<\mu$. As illustrated in \cref{introfig}(b), we would still have topological grounds to believe that population $x$ can remain bounded above zero in the long term, because $f(x_1)\pm \varepsilon>0$ and $f(x_2)\pm \varepsilon<0$. But under the same error tolerance $\epsilon$ on $\widehat f(y)$, we have $\widehat f(y)-\epsilon <0$ for all $y\in \R$ (\cref{introfig}(c)), so it is possible that no trajectories for population $y$ remain bounded above zero. In this way, the intensities $\mu$ and $\widehat\mu$ measure the degree to which  each attracting equlibrium ($x=1$ or $y=1$) persists to structure system behavior in the face of changes to the vector field---whether these are extrinsic, time-varying perturbations to population growth or intrinsic, time-independent changes to the model equations. \hfill //
\end{example}

In the body of this paper, we generalize the notion of intensity from \cref{ex:intro} to an analog for higher dimensions that still quantifies attractor persistence and is computable.
The organization is as follows.  

\Cref{sec:prelim} outlines the necessary preliminaries regarding $\R^n$ as a metric space and the dynamics generated by an autonomous differential equation $x'=f(x)$, $x\in\R^n$. 

In \cref{sec:control} we introduce time-dependent perturbations to the vector field $f$ via control functions $g(t)$. In this work control functions $g(t)$ play a role analogous to discrete errors in the theory of intensity of attraction for maps \cite{mcgehee1988}. In particular, studying the dynamics of $x'=f(x)+g(t)$ allows us to probe outwards from an attracting set that would otherwise be invariant under $f$ into the transient regions of its domain of atttraction.  We take controls $g$ to be essentially bounded, and endow them with an $L^\infty$ norm $\Infnorm{g}$ that reflects the norm in play on $\R^n$. In the sections that follow we quantify how large $\Infnorm{g}$ must become in order to escape the domain of attraction of an attractor under $f$---this is the essential ingredient for defining intensity of attraction.

Towards this end, \cref{sec:reach} broadens our attention from individual control functions to bounded collections of control functions. Given a set $S\subset\R^n$ and a bound $r$ on control, we ask which points in state space can be reached from $S$ using some control $g$ with $\Infnorm{g}<r$. The answer is called a reachable set, a term that carries a compatible and more general meaning in the theory of control \cite{colonius2000, maler} and differential inclusions \cite{BeynRieger2007}.  Depending on the vector field $f$ and the magnitude of $r$, the reachable set from an attractor might lie within its domain of attraction or spill outside this domain. 

\Cref{sec:IoA} defines intensity of attraction in terms of reachable sets:  roughly, intensity reflects the largest value of $r$ for which sets reachable from an attractor are bounded inside its domain of attraction. We connect intensity to resilience of predator-prey dyanmics in \cref{cycle}. Additional examples revisit \cref{introfig} with greater precision (\cref{indep}), illustrate how the value of intensity depends on the metric used to quantify perturbations (\cref{diagnorms}), and show how intensity relates to Colonius and Kliemann's point of first discontinuity in a map from control bounds to reachable sets \cite{colonius2000} (\cref{compareCK}).

\Cref{sec:continuation} demonstrates that intensity of attraction also quantifies attractor persistence in the face of time-\emph{independent} perturbations to a vector field.  \Cref{main} gives the main result that intensity bounds from below the distance in the space of vector fields that an attractor continues immediately. The key to this result is a link between reachable sets and the isolating blocks of Conley theory, which we develop in \cref{isolated} and \cref{reachblocks}. \Cref{ppcont} illustrates continuation of a predator-prey limit cycle through time-independent perturbations that don't exceed its intensity. We close by leveraging the present setting to provide a new proof that attractors are upper-semicontinuous (\cref{uppersc}).

\section{Preliminaries}
\label{sec:prelim}

\subsection{$\R^n$ as a metric space}
\label{subsec:metric}
Let $d:\R^n\times\R^n \rightarrow [0,\infty)$ be a homogeneous and translation-invariant metric on $\R^n$.
Let $\norm{\cdot}$ denote the corresponding norm with $\norm{x}=d(x,0)$ and $d(x,y)=\norm{x-y}$.  
Though we use the Euclidean metric and norm in most examples, we intentionally develop the theory for general $d$ and $\norm{\cdot}$.  In any particular application, context should inform the choice of metric and norm. 

The following definitions and notations are standard in this metric setting. We say a set $X\subset\R^n$ is bounded if there exists an $M>0$ such that $\norm{x}\leq M$ for all $x\in X$. 
The distance between two sets $X$, $Y\subset\R^n$ is \[\text{dist}(X,Y)\equiv\inf\{d(x,y) \ | \ x\in X, y\in Y  \}.\]  A single-point set $\{x\}$ is not distinguished from its element $x$.
For $\epsilon>0$ and $x\in\R^n$, we denote the $\epsilon$-ball about $x$ as \[B_{\epsilon}(x)\equiv\{y\in\R^n \ | \ d(x,y)<\epsilon\}.\]  
More generally, we denote the $\epsilon$-neighborhood of a set $S\subset\R^n$ as \[\mathcal{N}_\epsilon(S)\equiv\{x\in\R^n \ | \ d(x, S) < \epsilon\}.\]

The metric topology is used throughout. For $S\subset\R^n$, $\text{int}(S)$ and $\overline S$ denote the interior and closure of $S$, respectively. A set $N$ is a neighborhood of $S$ if $S\subset\text{int}(N)$. The following lemma states that a compact set has arbitrarily close compact neighborhoods; its straightforward proof is omitted.

\begin{lemma}\label{compactnbhd}
If $S\subset\R^n$ is compact, then for any neighborhood $N$ of $S$ there exists a compact set $K\subset \text{int}(N)$ that is also a neighborhood of $S$.
\end{lemma}

\subsection{Dynamics generated by $x'=f(x)$}
\label{subsec:dynamics}
Consider an autonomous system of ordinary differential equations 
\begin{equation}\label{f}
x'=f(x)
\end{equation}
where $x\in\R^n$, $x'$ denotes $\frac{dx}{dt}$, and $f$ maps from an open set $U\subset\R^n$ to $\R^n$.  
We assume that the vector field $f$ in the differential equation \cref{f} is bounded, so that $\fnorm{f}\equiv \sup\limits_{x\in U} \norm{f(x)}$ exists. We also take $f$ to be globally Lipschitz on $U$, so that there exists a unique local solution $x(t)$ in $U$ to each initial value problem \[x'=f(x), \ \ x(0)=x_0\in U.\] The solutions $x(t)$ generate a local flow $\varphi:\R\times U \supset D \to U$ given by $\flow{t}{x_0}=x(t)$. 
Though in general trajectories may leave $U$, preventing definition of a global flow $\varphi:\R\times U \to U$, we assume that $\varphi$ is defined on any time domain of interest.

Fixing the time coordinate of $\varphi$ yields a time-$t$ map $\varphi^t:U\rightarrow U$, given by $x\mapsto \flow{t}{x}$.  Functions such as the flow and the time-$t$ map are extended in the natural way to take sets as arguments. For example, if $T\subset\R$ and $S\subset\R^n$, then $\flow{T}{S}\equiv\bigcup\limits_{t\in T}\bigcup\limits_{x\in S} \flow{t}{x}$ and $\varphi^t(S)\equiv\bigcup\limits_{x\in S} \varphi^t(x)$.
We next put this notation to work in a few definitions that build to attractors, our main object of study.
\begin{definition}\label{inv}
 A set $S\subset U$ is \emph{forward invariant} under the flow $\varphi$ if $\varphi^t(S)\subset S$ for all $t>0$. $S$ is \emph{invariant} if $\varphi^t(S)= S$ for all $t\in\R$.
\end{definition}
\begin{definition}\label{omega}
The \emph{omega limit set} of a set $S\subset U$ is $\omega(S)\equiv \bigcap\limits_{T\geq 0}\overline{\bigcup\limits_{t\geq T} \varphi^t(S)}$.
\end{definition}
\begin{definition}\label{attr}
An \emph{attractor} is a compact, nonempty invariant set $A$ such that $A=\omega(N)$ for some neighborhood $N$ of $A$.
\end{definition}
The statement $A=\omega(N)$ indicates that the neighborhood $N$ is ``attracted'' asymptotically to $A$ in forward time. We call the collection of \emph{all} points that tend to $A$ in forward time $A$'s domain of attraction.
\begin{definition}\label{DofA}
The \emph{domain of attraction} of an attractor $A$ is \[\mathcal{D}(A)=\{x\in U \ | \  \varnothing\not=\omega(x)\subset A\}.\]
\end{definition}

\noindent Our primary objective is to quantify what we'll call the intensity of an attractor $A$, which should in some way reflect the strength of transient dynamics over $\domattr{A}$. To test that strength, we probe the domain of attraction with time-dependent control (\cref{sec:control}).

\section{Time-Dependent Control Functions}
\label{sec:control}

We introduce time-dependent perturbations to the system \cref{f} via nonautonomous control functions $g$, yielding a perturbed system
\begin{equation}
x'=f(x)+g(t). \label{eq:fg}
\end{equation}
The control functions $g$ are integrable and essentially bounded; that is, given a time interval $I\subset\R$ we take $g\in\gspaceI{I}$, equipped with the norm $\Infnorm{g}\equiv \inf \{ c\geq 0 \ : \ \mu\left(\{x\in I: \norm{g(x)}>c  \}\right)=0\}$. Note that the value of $\Infnorm{g}$, which plays an important role in \cref{sec:reach,sec:IoA,sec:continuation}, depends on the norm used in $\R^n$. 

With $g\in\gspaceI{I}$, the righthand side of \cref{eq:fg} satisfies the the Carath\'eodory conditions on $U\times I$: the function $F(x,t)=f(x)+g(t)$ is measurable in $t$ for fixed $x$, continuous in $x$ for fixed $t$, and (without loss of generality) bounded on any compact subset of $I\times U$ by the constant and hence integrable function $\norm{f}_\text{sup}+\Infnorm{g}$. Furthermore, $F$ inhertits $f$'s Lipschitz property in $x$. 
It follows (see, e.g. \cite{CL1955,Hale1980}) that for any initial condition $x(t_0) = x_0$ there exists a unique local solution $x(t)$ to \cref{eq:fg}, in the extended sense that 
\begin{align}
x(t)&=x(t_0)+\displaystyle\int_{t_0}^t \left( f(x(s))+g(s) \right) ds \label{eq:int} \\
\text{and }\hspace{1cm} x'(t) &= f(x) + g(t) \text{ almost everywhere.} \label{eq:ae}
\end{align}
The characterization \cref{eq:int} of solutions provides an analytic tool, while equation \cref{eq:ae} highlights the meaning of $\Infnorm{g}$ as a maximum deviation from the vector field $f$.  By default we assume that a unique local solution $x(t)$ to \cref{eq:fg} can be continued on any desired time interval $I$.

Existence and uniqueness of solutions to \cref{eq:fg} permit extension of standard flow notation to include the control function $g$ in a manner similar to \cite{colonius2000}:

\begin{definition}\label{gflow} For a given vector field $f$ and time interval $[0,T]$, the controlled flow 
\begin{align*}
\varphi:&[0,T]\times U \times \gspaceI{[0,T]} \rightarrow \R^n \\
\text{is given by } \hspace{2.2cm} &(\hat t,x_0,g)\mapsto \gflowf{\hat t}{x_0}{g}{f}=x(\hat t),
\end{align*}
 where $x(t)$ is the solution to $x'=f(x)+g(t)$, $x(0)=x_0$. 
\end{definition}
\noindent We omit the vector field $f$ when it is clear from context. To emphasize depedence on initial conditions, we write $\gflow{t}{\cdot}{g}=\gflowc{t}{\cdot}{g}$. In what follows, we move fluidly between describing solutions to $\cref{eq:fg}$ with the notation $x(t)$ and with the flow notation of Definition \ref{gflow}. The latter has the advantage of making $g$ visible. 

Continuity of $\gflow{t}{x}{g}$ with respect to $t$ and $x$ is well-known (see \cite{hsd2004}). 
Lemma \ref{contg} establishes continuity of $\gflow{t}{x}{g}$ with respect to control input $g$. 

\begin{lemma}\label{contg}
Fix a vector field $f:U\to\R^n$ and a time interval $[0,T]$, and let $g$ and $h$ be control functions in $\gspaceI{[0,T]}$. Then for any $\epsilon$, there exists a $\delta$ such that $\|g-h\|_\infty <\delta$ implies $\norm{\gflow{t}{x}{g}-\gflow{t}{x}{h}}<\epsilon$ for all $t\in[0,T]$ and all $x\in U$. 
\end{lemma}

\begin{pf}
Fix $\epsilon>0$. We apply Gr\"{o}nwall's Inequality to the function 
\[u(t)\equiv\norm{\gflow{t}{x}{g}-\gflow{t}{x}{h}}.\]
 Equation \cref{eq:int} gives that for any $x\in U$ and any $t\in[0,T]$,
\begin{align*}
\hspace {1cm} u(t)&=\norm{\int_0^t \big [ f(\gflow{s}{x}{g})+g(s) - f(\gflow{s}{x}{h}) -h(s) \big]  ds } \\
&\leq  \int_0^t \norm{ f(\gflow{s}{x}{g})- f(\gflow{s}{x}{h})  }  ds  + \int_0^t \norm{g(s) - h(s) }  ds \\
&\leq L \int_0^t u(s) ds  + T \|g-h\|_\infty
\end{align*}
where $L$ is the Lipschitz constant for $f$.
Gr\"{o}nwall's Inequality implies that \\ $u(t)\leq T\|g-h\|_\infty e^{Lt}$, and it follows that 
$\norm{\gflow{t}{x}{g}-\gflow{t}{x}{h}}\leq T\|g-h\|_\infty e^{LT}$.
Taking $\delta=\frac{\epsilon}{Te^{LT}}$ ensures that $\|g-h\|_\infty<\delta$ implies $\norm{\gflow{t}{x}{g}-\gflow{t}{x}{h}}<\epsilon$ for all $t\in[0,T]$ and all $x\in\R^n$, as desired. \qed
\end{pf}

\section{Reachable Sets}
\label{sec:reach}

The controlled flow of \cref{gflow} gives outcomes under a single function $g\in\gspaceI{I}$, which highlights the effect of a known intervention in a controlled system. However, when $g$ represents an uncertain disturbance, it is natural to consider trajectories corresponding to an entire family of such functions. We restrict our attention to families of the form
\[B_r(0)\equiv \{g\in\gspaceI{I} \ | \ \Infnorm{g}\leq r\}. \]

\noindent By using the collection $B_r(0)$, we focus on situations in which bounds on disturbances are known, even though their exact forms are not. This approach is compatible with studies of random dynamical systems  with bounded noise (e.g. \cite{homburg2010bifurcations, lamb2015topological, zmarrou2007bifurcations}). But in contrast to systems with unbounded noise such as stochastic differential equations, trajectories cannot necessarily reach any state from another.

We use $\reach{r;f}{S}$ to denote the collection of states that can be reached in forward time from a point in $S\subset U$ by modifying the vector field $f$ with $r$-bounded control:

\begin{definition}\label{PrS}
$\reach{r;f}{S}=
\bigcup\limits_{t\geq 0} \bigcup\limits_{x\in S} \ \bigcup\limits_{\Infnorm{g}\leq r} \gflowf{t}{x}{g}{f}$.
\end{definition}

\noindent When no confusion could result by omitting the vector field $f$, we write $\reach{r}{S}$.
In the present work we'll consider sets reachable from an attractor $A$. These can be found easily for one-dimensional systems, as the following proposition shows.

\begin{proposition}\label{PrAmedium}
Consider the one-dimensional system $x'=f(x)$ with $f(0)=0$, $f>0$ on $(-\epsilon,0)$, and $f<0$ on $(0,\epsilon)$ for some $\epsilon>0$. The reachable set from the attractor $A=\{0\}$ is $\reach{r}{A}=(a,b)$ where 
\begin{align*}
a&=\begin{cases}\sup\{x<0 \ | \ f(x)=r\} & \text{if } f(x)=r \text{ for some } x<0 \\ -\infty & \text{otherwise} \end{cases} \\
\text{and }\hspace{1em} b&=\begin{cases}\inf\{x>0 \ | \ f(x) = -r\}& \text{if }f(x)=-r \text{ for some } x>0 \\ \infty &\text{otherwise}\end{cases}.
\end{align*} 
\end{proposition} 
\begin{pf}
We begin with the inclusion $(a,b)\subset \reach{r}{A}$. Let $c$ be any point in $[0,b)$ and let $-m$ be the minimum value of $f$ on $[0,c]$. This minimum must be strictly greater than $-r$ by construction of $b$. Let $g$ be the constant function $g(t)=(m+r)/2$, so that $0<m<g(t)<r$. Then $\|g(t)\|_\infty<r$. Also, since $f(x)\geq -m$ on $[0,c]$, $g(t)>m$ implies $f(x)+g(t)>0$ on $[0,c]$. Using separation of variables, one may confirm that $c=\gflowc{t}{0}{g}$ for the finite, non-negative time $t=\int_0^c \frac{dx}{f(x)+(m+r)/2}$. A similar argument gives that any point in $(a,0]$ is reachable in finite time by a constant control strictly bounded by $r$. Definition \ref{PrS} then implies that $(a,b)\subset \reach{r}{A}$. 

The reverse inclusion follows trivially if $a=-\infty$ and $b=\infty$, but requires more work than one might expect in the case that $a$ or $b$ is finite.  Suppose $b<\infty$. Since $f$ is continuous, $f(b)=-r$ and $f(b)+r=0$. Therefore, choosing $g(t)\equiv r$ in equation \cref{eq:fg} yields an autonomous flow with $b$ as an equilibrium.

The Lipschitz property of $f$ implies that $f(x)+r\leq -L(x-b)$ on $[0,b]$ for some $L>0$.
We use this fact and a comparison of solutions to the initial value problems
\begin{equation}\label{LipDE}
x'=-L(x-b), \hspace{0.2cm} x(0)=0
\end{equation} and 
\begin{equation}\label{gDE}
y'=f(y)+g(t), \hspace{0.2cm} y(0)=0
\end{equation}
to show that $\gflowc{t}{0}{g} < b$ when  $\Infnorm{g}\leq r$.
For brevity, let $x(t)$ denote the solution to \cref{LipDE} and $y(t)$ denote the solution to \cref{gDE} for a fixed $g$.

Assume for the sake of contradiction that there exists a $t>0$ and $g\in\gspaceI{[0,t]}$, $\|g\|_\infty\leq r$, such that %$\gflow{t}{0}{g}\geq b$; i.e. 
$y(t)\geq b$. Let $T$ denote the minimum positive time $t$ at which $y(t)=b$. 
Since $x(t)$ remains strictly bounded below $b$ and both paths are continuous, the inequality $y(t)\geq x(t)$ must hold on some time interval $[t_*,T]$ with $0\leq t_*<T$. In fact, one may choose $t_*$ so that $x(t_*)=y(t_*)$. For $s\in[t_*,T]$, these inequalities must then hold: 
\begin{equation}\label{longIneq}
f(y(s))+g(s) \leq f(y(s))+r \leq -L(y(s)-b) \leq  -L(x(s)-b)
\end{equation}
The second inequality follows from the Lipschitz property of $f(\cdot)+r$, which gives that $f(y(s))+r
\leq|(f(y(s))+r)-(f(b)+r)| \leq L |y(s)-b|=-L(y(s)-b)$ for $s\in[t_*,T]$.\\

On the other hand, the assumptions $x(T)<y(T)$ and $x(t_*)=y(t_*)$ imply that 
\begin{equation}
\int_{t_*}^T-L(x(s)-b)ds <\int_{t_*}^T\left[f(y(s))+g(s)  \right]ds \label{intIneq}
\end{equation}
Comparing inequalities \cref{longIneq} and  \cref{intIneq} yields the desired contradiction. Hence the solution $y(t)$ to IVP \cref{gDE} must remain strictly bounded below $b$ in forward time; i.e. $\gflowc{t}{0}{g}\in (-\infty,b)$ for $\Infnorm{g}\leq r$ and $t\geq0$. \Cref{PrS} implies that $\reach{r}{A} \subset (-\infty, b)$. 
In the case that $a>-\infty$, the inclusion $\reach{r}{A}\subset (a, \infty)$ follows similarly. Together, the inclusions $\reach{r}{A} \subset (-\infty, b)$ and $\reach{r}{A}\subset (a, \infty)$ imply that $\reach{r}{A}\subset (a,b)$ when $a$, $b$, or both are finite. \qed
\end{pf}

In the case of a one-dimensional linear system, \cref{PrAmedium} implies that the set reachable from the origin scales in direct proportion to the control bound and in inverse proportion to the eigenvalue magnitude:

\begin{corollary}\label{PrAeasy}
For the one-dimensional, linear system $x'=-\lambda x$ ($\lambda>0)$ with global attractor $A=0$, the reachable set from $A$ for $r>0$ is $\reach{r}{A}=(-\frac{r}{\lambda} ,\frac{r}{\lambda} )$. 
\end{corollary}

Reachable sets are, in general, harder to compute than in \cref{PrAmedium} and \cref{PrAeasy}. The following example shows that the computation can be non-trivial even for a diagonal, two-dimensional linear system.

\begin{example}\label{PrAhard} The system 
\begin{align}
x' &= -x\\
y' &= -2y
\end{align}
on $\R^2$ has global attractor $A=(0,0)$. Considering each variable separately, \cref{PrAeasy} indicates that $\reach{1}{A}$ is contained in the rectangular region $R=(-1,1)\times(-\frac{1}{2},\frac{1}{2})$. Furthermore, the family of constant control functions $\{g_{\theta,c}(t)=c\cos\theta+c\sin\theta \ | \ \theta\in[0,2\pi), 0\leq c \leq 1\}$ sends the origin asymptotically to points on ellipses $x^2+4y^2=c^2$. Therefore, the reachable set $\reach{1}{A}$ contains the open region $\mathcal{O}$ enclosed by the ellipse $E=\{(x,y)\ | \ x^2+4y^2 =1 \}$. Control functions with magnitude $1$ directed opposite to the vector field also drive trajectories asympototically from $\mathcal{O}$ to $E$. 

\begin{figure}[h]
\centering
\includegraphics[scale=0.9]{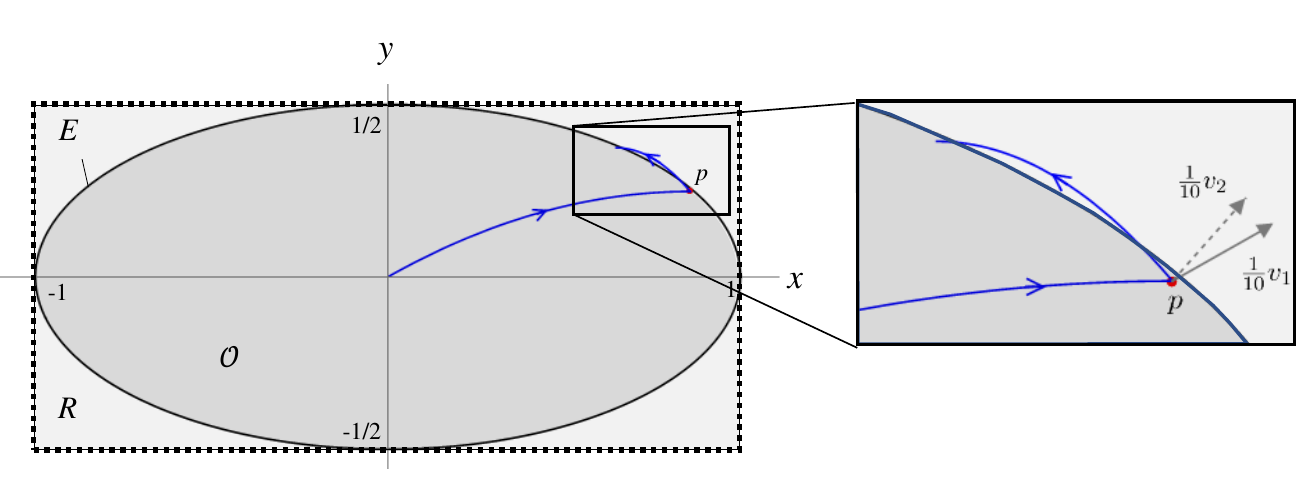}
\caption{Inexact bounds on a 2D reachable set. The set in question consists of points reachable from the origin in the linear system $x'=-x$, $y'=-2y$ with control bound $r=1$. This set must be nested between the open set $\mathcal{O}$ and the rectangle $R$. Piecewise-constant control that shifts directions from $v_1$ to $v_2$ drives a trajectory (blue) outside of $\mathcal{O}$.}

\label{2DlinearPrA}
\end{figure}

Because the vector field has magnitude $1$ on $E$, one might expect that $\reach{1}{A}=\mathcal{O}$. But the vector field is not normal to $E$ except at the vertices of the ellipse. Therefore, adjusting a control function's direction to push perpendicular to $E$ can drive a solution outside $\mathcal{O}$. Figure \ref{2DlinearPrA} illustrates this phenomenon for piecewise constant control 
\begin{equation}
g(t)=\begin{cases} v_1& \text{ if } 0\leq t<10\\ v_2 & \text{ if } 10\leq t\end{cases}
\end{equation}
where $v_1=0.99 \left[ \begin{array}{c} \sqrt{3}/2 \\ 1/2 \end{array}\right]$ and $v_2=0.99\left[\begin{array}{c} \sqrt{3/7}\\ 2/\sqrt{7}\end{array}\right] $  have magnitude 0.99 but different directions.
We therefore have the bounds $\mathcal{O}\subsetneq\reach{1}{A}\subset R$. Athough these bounds follow directly from analytic computation of trajectories and \cref{PrAeasy}, a closed form for $\reach{1}{A}$ does not.  \hfill //
\end{example}

Numerical algorithms for computing reachable sets are an active topic of research in the fields of control systems, differential inclusions, and optimization. The introductory piece \cite{maler} outlines fundamental discretization challenges nicely.  Existing approaches  include set-valued Euler schemes \cite{BeynRieger2007,Rieger2015, Rieger2016, Sandberg2008}, optimal control algorithms \cite{Rasmussen2017,Baier2013}, and level-set methods (\cite{deWeerdt2011} and references therein).

Though the approximation of reachable sets can be involved, analytically they obey nice nesting properties. The following lemma, which follows immediately from \cref{PrS}, states that the real number ordering of control bounds induces an order by inclusion on reachable sets:

\begin{lemma}\label{rnest2}
If $r < r'$ then $\reach{r}{S}\subset \reach{r'}{S}$.
\end{lemma}

\noindent A stronger statement holds when the set $S$ is an attractor $A$:
\begin{proposition}\label{PrAtoA}
 For $A$ an attractor, $\reach{r}{A}\searrow A$ as $r\searrow 0$, in the sense that 
 \begin{enumerate}
 \item[(1)] $r < r'$ implies $\reach{r}{A} \subset \reach{r'}{A}$
 \item[(2)] for any neighborhood $V$ of $A$, there is an $r>0$ such that $\reach{r}{A}\subset V$, and
 \item[(3)] $A=\bigcap_{r> 0} \reach{r}{A}$.
 \end{enumerate}
 \end{proposition}
 
\noindent We delay the proof of \cref{PrAtoA} until \cref{upper}, where it follows from existence of attractor blocks introduced in \cref{isolated} and supports a new proof of the semicontinuity of attractors.
 
\section{Intensity of Attraction}
\label{sec:IoA}
Intensity of attraction, introduced by McGehee for maps \cite{mcgehee1988}, carries over naturally to the flow setting by replacing sets reachable from an attractor under $\epsilon$ pseudo-orbits with sets reachable from an attractor under $r$-bounded control.

\begin{definition}\label{Idef}
The \emph{intensity of attraction} of an attractor $A$ is 
\begin{equation}
\intensity{A}=\sup \{ r \geq 0 : P_r(A)\subset K \subset \mathcal{D}(A) \text{ for some compact } K\subset\R^n\}
\end{equation}
\end{definition}

\noindent In other words, intensity of attraction reflects the control magnitude $\Infnorm{g}$ necessary to escape from every compact subset of an attractor's domain. 

\begin{proposition}\label{mug0}
For $A$ an attractor, $\intensity{A}>0$.
\end{proposition}

\noindent Our proof of \cref{mug0} relies on attractor blocks and \cref{PrAtoA} and appears in \cref{upper}.

The following two examples, though somewhat synthetic, serve to illustrate two points about attractor intensity: it is independent of eigenvalues and basin size (\cref{indep}), and its value depends on the norm $\norm{\cdot}$ in play on $\R^n$ (\cref{diagnorms}).

\begin{example}\label{indep}
It is possible for two vector fields that share the same attracting equilibrium $A$, the same eigenvalues at $A$, and the same domain of attraction $\domattr{A}$ to nonetheless exhibit different intensities of attraction $\mu(A)$. Consider the pair of $C^1(\R,\R)$ vector fields plotted in \cref{introfig}:
\begin{align}
f(x)&=x-x^2\\
\text{and } \ \hat f(x)&=\begin{cases} x & \text{if } x<0\\ \pi^{-1}\sin(\pi x) & \text{if } 0\leq x < 1 \\ 1-x &\text{if } x\geq1 \end{cases}
\end{align}
The equilibrium $1$ is an attractor for both $f$ and $\hat f$, with eigenvalue $-1$. Its domain of attraction is $(0,\infty)$ in each case. In fact, the eigenvalue at the repelling equilibrium $0$ that consitutes the boundary of $\domattr{A}$ is also $1$ for both $f$ and $\hat f$. Yet $\mu(A)=1/4$ under $f$ while $\hat\mu(A)=\pi^{-1}\approx0.318$ under $\hat f$. This example alerts us that eigenvalues and size of a domain of attraction---two common measures of resilience in engineering and ecology \cite{meyer2016}---can neglect differences in the strength of transient dynamics over $\domattr{A}$---a feature that intensity of attraction detects. \hfill //
\end{example}

\begin{example}\label{diagnorms}
To easily see how $\mu(A)$ depends on norm choice, consider the system
\begin{align}
x'&=\frac{\sqrt{2}}{4}(x+y)^2-x \label{x}\\
y'&=\frac{\sqrt{2}}{4}(x+y)^2-y \label{y}
\end{align}
whose phase portrait is given in \cref{diag}(a). Rotating coordinates by $\pi/4$ to $u=\frac{\sqrt{2}}{2}(x+y)$ and $v=\frac{\sqrt{2}}{2}(y-x)$ gives \begin{align}
u'&=u(u-1) \label{u}\\
v'&=-v.  \label{v}
\end{align}
In these coordinates it is clear that $A=(0,0)$ is an attractor with domain of attraction $\mathcal{D}(A)=\{(u,v) \ | \ u <1\}$. The uncoupled form of equations \cref{u} and \cref{v} and the geometry of $\mathcal{D}(A)$ allow us to restrict attention to the invariant line $v=0$ when calculating $\mu(A)$. On the interval $0\leq u \leq 1$, the  maximum strength of the vector field in the direction of the origin is $1/4$, achieved by the vector $m$ at $u=1/2$. A control function must push with magnitude greater than $1/4$ in the direction of $-m$ to steer a trajectory from $A$ to $\mathcal{D}(A)^C$. The critical control bound needed for this push, $\mu(A)$, depends on the norm on $(x,y)$ space. Under the Euclidean norm (the 2-norm), the intensity is $\mu_2(A)=1/4$ 
(\cref{diag}(c)). 
The intensity of $A$ under the 1-norm is $\mu_1(A)=\sqrt{2}/4$ (\cref{diag}(b)) and the intensity under the max norm is $\mu_\text{max}(A)=1/4\sqrt{2}$ (\cref{diag}(d)). More generally, a straightforward calculation gives that under the $p$-norm $\norm{(x,y)}_p=(x^p+y^p)^{1/p}$, the intensity of $A$ is $\mu_p(A)=2^{(\frac{1}{p}-\frac{1}{2})}\cdot\frac{1}{4}$. \hfill //

\begin{figure}[h!]
\centering
\includegraphics[scale=0.6]{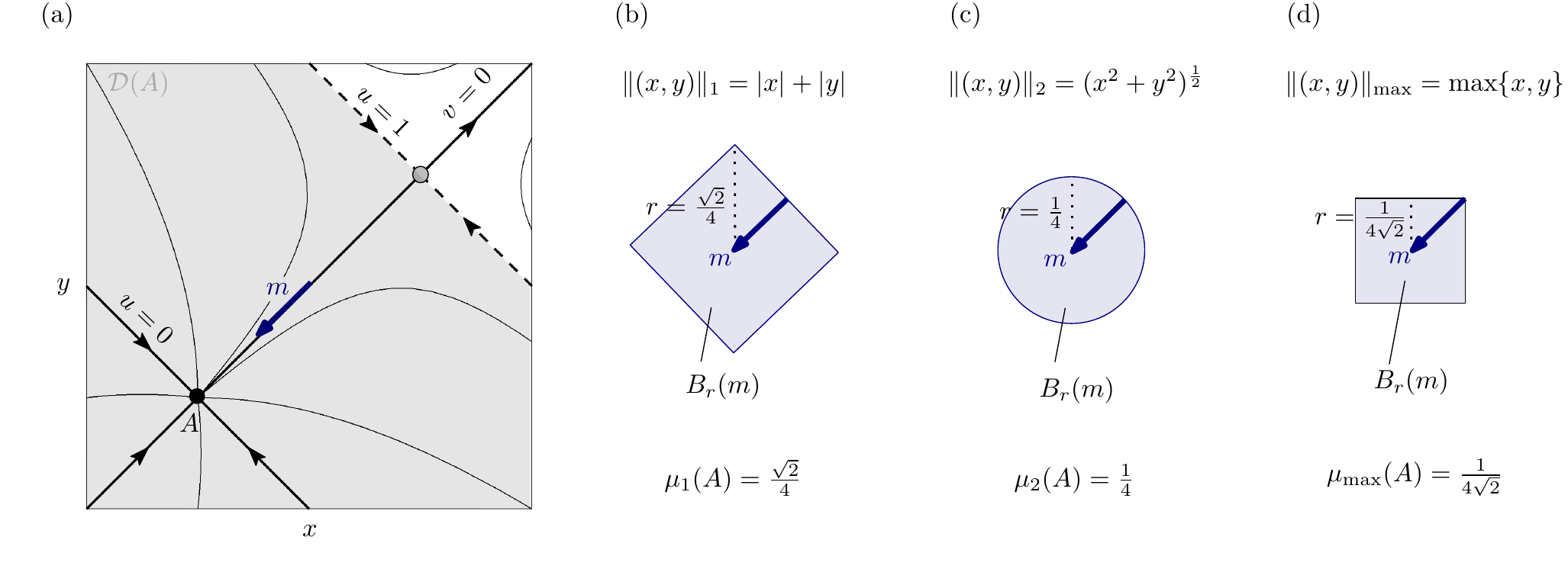}
\vspace{-0.5cm}
\caption{Intensity depends on norm / metric. (a) Phase portrait for the system \cref{x,y}. $A$ is an attracting equilibrium and the dashed line marks the boundary of its domain of attraction, shaded grey. Analysis in rotated coordinates (equations \cref{u,v}) reveals that the intensity of $A$ is the norm of the vector $m$. (b), (c), (d): The 1-norm, 2-norm, and max norm each yield different values of intensity of attraction for $A$.}
\label{diag}
\end{figure}
\end{example}

In specific applications, one should choose a norm compatible with the relevant notions of disturbance or uncertainty.  We anticipate that the max norm will be appropriate when bounds on vector field perturbations are known componentwise rather than jointly. As in \cref{diagnorms}, one expects $\mu_{p}(A) \leq \mu_{q}(A)$ when $p>q$. Therefore, choosing a max norm may also give the most conservative estimate of attractor intensity.   

We next compare intensity to a closely related concept examined by Colonius and Kliemann \cite{colonius2000}: the lowest point of discontinuity in a map from control bound to the reachable set (control set) of an attractor. \cref{compareCK} demonstrates that this first discontinuity may differ from the intensity of the attractor. 

\begin{example}\label{compareCK} 
The one-dimensional system with vector field $f(x)=\frac{3}{4}x^4-x^3-3x^2-1$ has an attractor $A=x_0$ and repeller $R=x_1$
 corresponding to the roots $x_0$ and $x_1$ of the quartic (\cref{compareCKfig}(a)). The domain of attraction of $A$ is $\mathcal{D}(A)=(-\infty,x_1)$. 
 Consider the map $C:[0,\infty)\rightarrow
 \mathscr{P}(\R)$
 given by $r\mapsto \reach{r}{A}$. \cref{compareCKfig}(b) depicts the graph of $C$.
As the control bound $r$ increases from zero, the first discontinuity in $C$ occurs at $r=2.25$, corresponding to the local minimum of $f$ at $(-1,-2.25)$. However, reachable sets remain bounded within $\mathcal{D}(A)$ for small enough $r>2.25$. It is not until the second discontinuity at $r=9$, corresponding to the global minimum of $f$ at $(2,-9)$, that reachable sets escape $\mathcal{D}(A)$. Hence the intensity $\mu(A)=9$ exceeds the value of $r$ at the first discontinuity ($\rho^*=2.25$ in the notation of \cite{colonius2000}). \hfill //

\begin{figure}[h]
\centering
\includegraphics[scale=0.8]{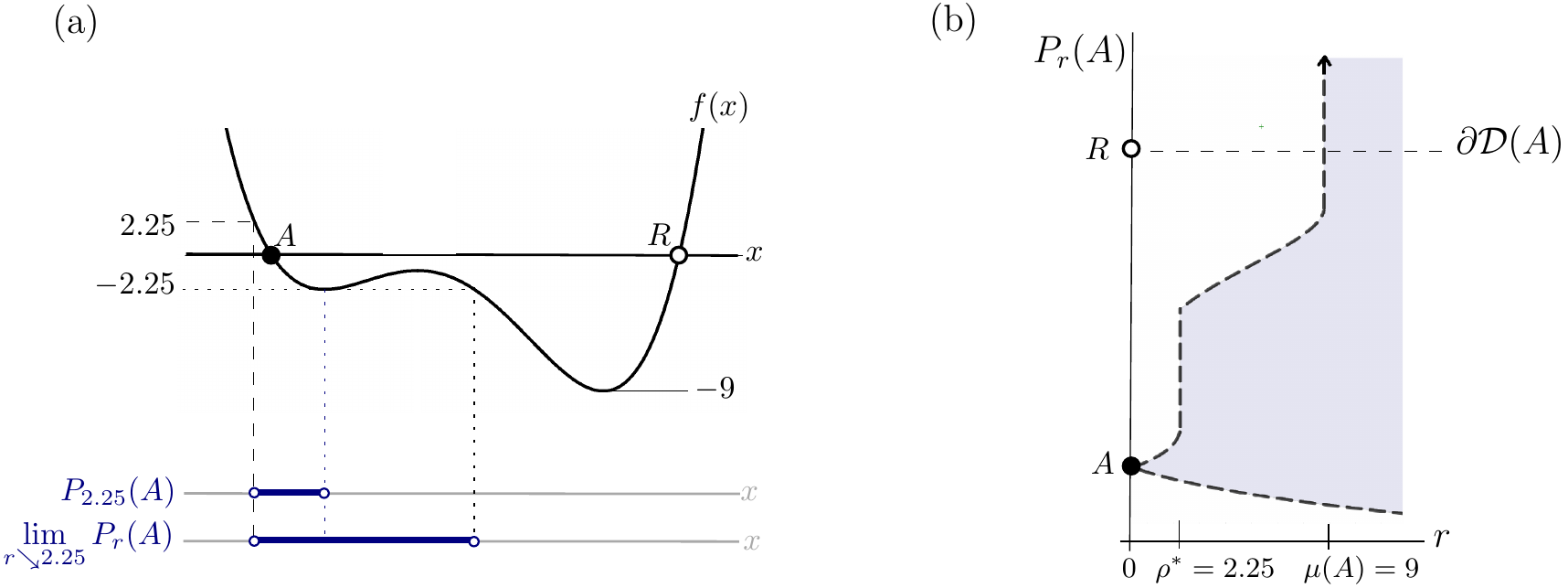}
\caption{Intensity versus first discontinuity.
 (a) The vector field $f$ of \cref{compareCK} yields an attractor $A$ and repeller $R$. As the control bound $r$ increases, the first discontinuity in the map $C$ that sends $r$ to $P_r(A)$ occurs at $\rho^*=2.25$.  (b) The graph of $C$ illustrates that the first discontinuity at $\rho^*=2.25$ occurs well below the intensity $\mu(A)=9$. }
\label{compareCKfig}
\end{figure}

\end{example}

\Cref{compareCK} shows that it is possible for reachable sets to expand discontinuously without escaping a domain of attraction. Whether the reverse is possible---reachable sets escaping a domain of attraction without expanding discontinuously---is an interesting question not pursued further here.

We close this section with an ecological application that illustrates attractor intensity for a stable limit cycle, and the information about resilience it carries. 
\begin{example}\label{cycle}
Consider the following model of predator ($y$) and prey ($x$) dynamics, the details of which can be found in \cite{may1972,rosenzweig1971}:
\begin{align}
\frac{dx}{dt}&=ax\left(1-\frac{x}{K}\right)-ky(1-e^{-cx}) \label{pp1}\\
\frac{dy}{dt}&=-by+\beta y(1-e^{-fx}). \label{pp2}
\end{align}
 
As the prey carrying capacity $K$ grows from $3$ to $4$, a Hopf bifurcation in the first quadrant transforms the stable equilibrium $A_1$  (\cref{cyclefig}(a)) into an unstable equilibrium and a stable limit cycle $A_2$ (\cref{cyclefig}(c)).  Rosenzweig used this predator-prey system to illustrate a ``paradox of enrichment''  in which adding nutrients to an ecosystem (thus increasing $K$) destabilizes the ecosystem \cite{rosenzweig1971}. His analysis was restricted to the dynamics near the rest point in the first quadrant, which indeed turns from stable to unstable as $K$ increases. But in what sense is the \emph{periodic orbit} $A_2$ in Figure \ref{cyclefig}c ``destabilized'' relative to the spiral sink $A_1$ in \cref{cyclefig}(a)? Both invariant sets are stable by standard mathematical definitions. In present vocabulary we might say adding nutrients lowered the resilience of the populations. Intensity of attraction provides one measure of this change in resilience. 

\Cref{cyclefig} plots the result of reachable set computations for both the spiral sink (panel b) and stable limit cycle (panel d). Reachable sets were computed using a fixed-grid set-valued Euler
\begin{figure}[h!]
\centering
\includegraphics[width=0.76\textwidth]{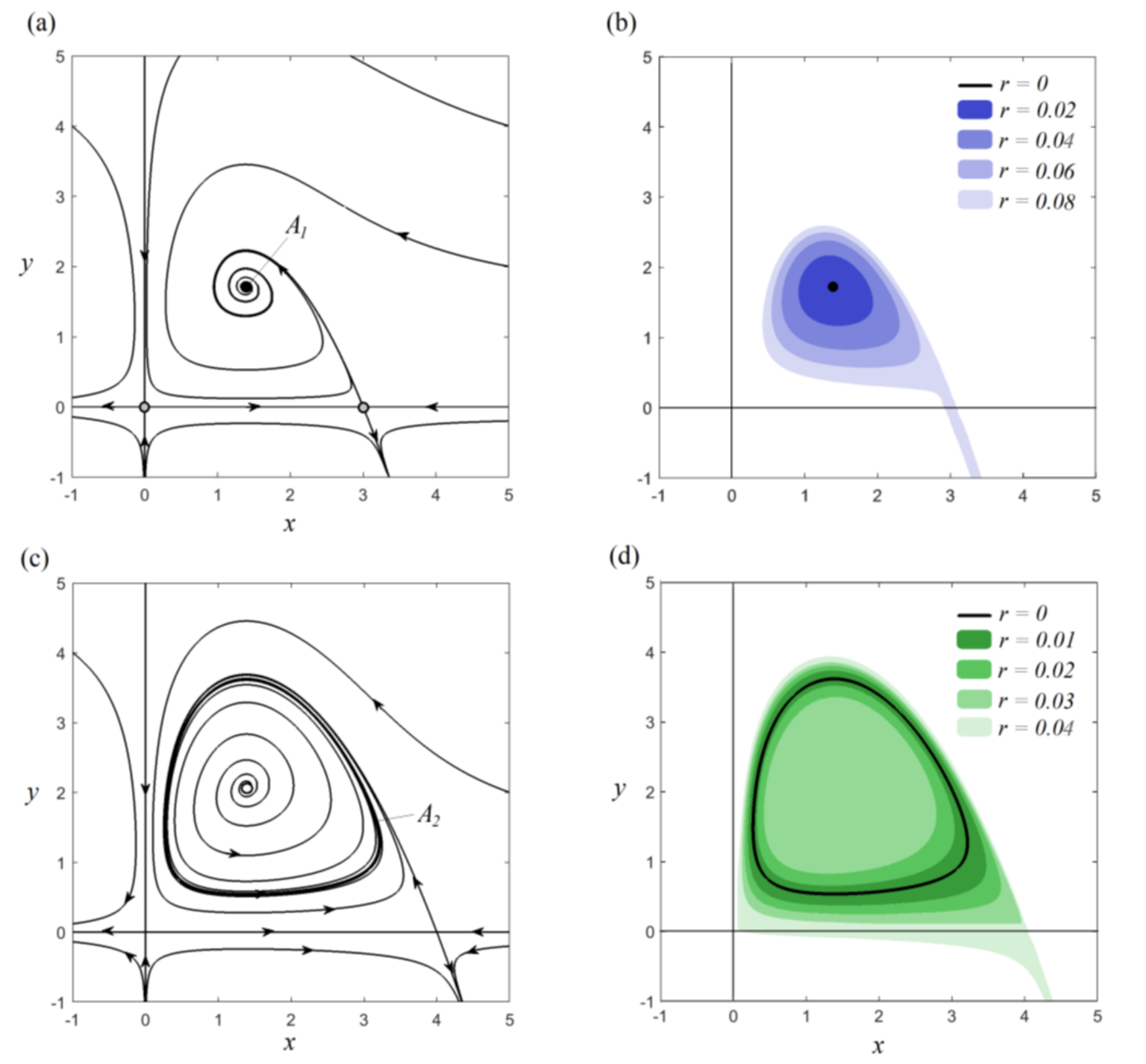}
\caption{Phase portraits (a,c) and reachable sets (b,d) for a predator-prey equilibrium (a,b) and limit cycle (c,d) corresponding to equations \cref{pp1,pp2}.  Parameters, in arbitrary units, are $K=3$ (a, b), $K=4$ (c,d), and $a=1$, $k=0.5$, $c=1.5$,  $b=0.5$, $\beta=1$, $f=0.5$ (all panels).}
\label{cyclefig}
\end{figure}

\noindent method at a spatial resolution of $10^{-3}$. Sets reachable from the spiral sink remain bounded within its domain of attraction (the first quadrant) for $r\leq 0.06$. For $r=0.08$, the reachable set spills over the $x$-axis. From \cref{cyclefig}(b), we estimate $0.6<\mu(A_1)<0.8$. 
\cref{cyclefig}(d) illustrates sets reachable from the periodic orbit $A_2$. These remain bounded within its domain of attraction (the first quadrant minus the unstable rest point) for $r\leq 0.02$; for $r=0.03$, the reachable set includes the unstable rest point. Hence $0.02<\mu(A_2)<0.03$. If we instead consider the attractor $A_3$, consisting of the periodic orbit $A_2$ and the region it encloses, we get a slightly higher intensity $0.03<\mu(A_3)<0.04$. Regardless of whether we consider $A_2$ or $A_3$, the intensity estimates agree with Rosenzweig's argument that nutrient enrichment lowers stability (resilience): the attractor in the first quadrant has a lower intensity of attraction in the ``enriched'' scenario $K=4$ compared to $K=3$. As a result, we might expect that environmental variability is more likely to drive the predator and/or prey population to extinction from the cycle $A_2$ than from the equilibrium $A_1$. \hfill //
\end{example}

\section{Continuation of Attractors}
\label{sec:continuation}

We now turn from time-dependent to time-independent perturbations of a vector field. \Cref{isolated} introduces basic definitions and results regarding isolated invariant sets in general and attractors in particular. \Cref{reachblocks} establishes a connection between reachable sets and attractor blocks, which we exploit in \cref{dist} to show that intensity of attraction gives a lower bound on attractor continuation distance in the space of vector fields. \Cref{main}, our main result, grounds the abstract question of attractor persistence in infinite dimensional vector field space to the study of a single system's reachable sets. We conclude in \cref{upper} with a proof of the semicontinuity of attractors based on reachable sets.

\subsection{Attractors as Isolated Invariant Sets}
\label{isolated}

Attractors are a special type of isolated invariant set, a useful object introduced by Conley for studying robust qualitative features of flows \cite{conley1978,conleyeaston1971}. The theory of isolated invariant sets has been developed in the settings of flows on smooth manifolds \cite{conleyeaston1971,wilsonyorke1973} and locally compact metric spaces \cite{mischaikowmrozek1999}. Definitions \ref{in} and \ref{iis} tailor those of \cite{mischaikowmrozek1999} to flows on $\R^n$.

\begin{definition}\label{in}
An \emph{isolating neighborhood} for a flow $\varphi$ is a compact set $N\subset\R^n$ whose invariant part $\text{Inv}(N,\varphi)\equiv\{x\in N \ | \ \flow{\R}{x}\subset N \text{ for all } t\in\R\}$ lies in the interior of $N$.
\end{definition}
\begin{definition}\label{iis}
A set $I\subset\R^n$ is an \emph{isolated invariant set} if $I=\text{Inv}(N)$ for some isolating neighborhood $N$.
\end{definition}

\begin{remark}\label{aiis}
\Cref{attr,iis} imply that an attractor $A$ is an isolated invariant set; the isolating neighborhood $N$ can be constructed via Lemma \ref{compactnbhd} as a compact set nested between $A$ and $\hat{N}$, where $\hat{N}$ is a neighborhood of $A$ such that $A=\omega(\hat{N})$. 
\end{remark}

An important consequence of \cref{in} is that isolating neighborhoods continue as such under sufficiently small perturbations to the flow, though the maximal invariant set in the interior may change \cite{mischaikowmrozek1999}. Homotopic and homological Conley indices built from an isolating neighborhood can yield coarse information about the isolated invariant set inside. These allow one to deduce topological features of an isolated invariant set that persist as it continues in nearby systems.

When the flow interacts nicely with the boundary of an isolating neighborhood, that neighborhood earns the name \emph{isolating block}. Various definitions of ``nicely'' have been used in \cite{conley1978,conleyeaston1971,mischaikowmrozek1999,wilsonyorke1973}. Here we restrict our attention to isolating blocks associated with attractors, henceforth called attractor blocks. Roughly, the vector field should point inward to the interior of an attractor block at each point along its boundary. The following definition generalizes this notion slightly.
\begin{definition}
An \emph{attractor block} $B$ for a flow $\varphi:\R\times\R^n\to\R^n$ is a nonempty, compact subset of $\R^n$ satisfying $\varphi^t(B)\subset \text{int}(B)$ for all $t>0$.
\end{definition}

\noindent Any attractor block is associated with an attractor in its interior:

\begin{lemma}\label{obattr}
For $B$ an attractor block, $\omega(B)$ is an attractor.
\end{lemma}

\begin{pf}
$\omega(B)$ must lie interior to $B$ because $\varphi^t(B)\subset\text{int}(B)$ for all $t>0$. Therefore $\omega(B)$ is the omega limit set of its neighborhood $B$. Further, it follows directly from \cref{omega} that $\omega(B)$ is invariant and closed; as a subset of compact $B$ it must also be compact. By \cref{attr}, $\omega(B)$ is an attractor. \qed
\end{pf}

\begin{definition}
If $B$ is an attractor block and $\omega(B)=A$, we call $B$ an attractor block \emph{associated with} $A$.
\end{definition}

\noindent We have already seen in \cref{obattr} that any attractor block is associated with an attractor. \Cref{ab} states that the converse also holds.
\begin{lemma}\label{ab}
For any attractor $A$ and any neighborhood $V$ of $A$, there exists an attractor block $B\subset \text{int}(V)$ associated with $A$.
\end{lemma}

We sketch a construction of such an attractor block using a Lyapunov function. Fix a neighborhood $V$ of $A$. Let $N$ be another neighborhood of $A$ such that $A=\omega(N)$ and let $\hat{N}$ be an isolating neighborhood with $\text{Inv}(\hat{N})=A$ and $\hat{N}\subset N$ (see Remark \ref{aiis}). Wilson and Yorke have shown there is an open neighborhood $\Omega$ of $A$ in $\hat{N}$ and a smooth, monotone Lyapunov function $L:\Omega\rightarrow \R$ such that $L(x)=0$ if $x\in A$ and $\frac{d}{dt}L(x(t))<0$ if $x\in\Omega - A$ (\cite{wilsonyorke1973}, Theorem 2.2). By Lemma \ref{compactnbhd} there exists a compact neighborhood $W$ of $A$ that lies inside $\text{int}(V\cap \Omega)$. Let $\epsilon=\min_{x\in\partial W}\{L(x) \}$. 
Let $B=L^{-1}([0,\epsilon])$. Strict monotonicity of the Lyapunov function implies both that $B$ is an attractor block ($\varphi^t(B)\subset\text{int}(B)$ for all $t>0$) 
and that $B\subset W\subset\text{int}(V\cap \Omega)$. 
It follows that $B\subset \text{int}(V)$. Finally, the attractor block $B$ is associated with $A$ because $A\subset B \subset N$, implying $\omega(B)=A$.

\subsection{Persistent Attractor Blocks from Reachable Sets}
\label{reachblocks}
In this subsection we show that one can construct attractor blocks from certain reachable sets; furthermore, these blocks persist as such not just for sufficiently small perturbations to the vector field, but for perturbations smaller than $r$, the metric bound on control.  Before stating these results we establish a definition and technical lemma.

Fixing all but the third argument in the controlled flow function $\gflowf{t}{x}{g}{f}$ of \cref{gflow} yields a map from control functions to trajectory endpoints:
\begin{definition}\label{G}
Let $\cntrlmapf{T}{x_0}{f}:\gspaceI{[0,T]} \to \R^n$ be given by $g\mapsto \gflowf{T}{x_0}{g}{f}$. The subscripts on $G$ may be omitted when $f$, $T$, and $x_0$ are clear from context, or general.
\end{definition}

\begin{lemma}\label{openG}
$\cntrlmapp$ is an open map. 
\end{lemma}
\begin{pf}
It suffices to show that $\cntrlmapp$ maps any basis element in the metric topology on $\gspaceI{[0,T]}$---an open ball $B_r(g)\equiv \{h\in\gspaceI{[0,T]} : \| h - g \|_\infty < r \}$---to an open set in $\R^n$. 
To show that $\cntrlmapp(B_r(g))$ is open, we will show that for any $h\in B_r(g)$ there exists an $\epsilon_h$ sufficiently small so that the open $\epsilon_h$-ball around $\cntrlmapp(h)$ is also in the image $\cntrlmapp(B_r(g))$. 

Fix $h\in B_r(g)$ and consider a point $\cntrlmapp(h) + v$ in $\R^n$. We derive an $\epsilon_h$ so that $\norm{v}<\epsilon_h$ implies there exists $k\in B_r(g)$ such that $\cntrlmapp(k)=\cntrlmapp(h)+v$. First construct a path from $x_0$ to $\cntrlmapp(h)+v$: let $x_h(t)=\gflow{t}{x_0}{h}$ and let
\[ \tilde x(t)=x_h(t)+t\frac{v}{T}. \]

\noindent Then $\tilde x(0) = x_h(0)=x_0$ and $\tilde x(T)=x_h(T)+v=\cntrlmapp(h)+v$, so $\tilde x$ is a path from $x_0$ to $\cntrlmapp(h)+v$. The velocity along the path $\tilde x$ (defined almost everywhere) is
\[\frac{d}{dt}\tilde x(t) = \frac{d}{dt}\left[ x_h(t) \right] +\frac{v}{T}=f(x_h(t))+h(t) + \frac{v}{T}.\]

\noindent The control $k(t)$ required to achieve the path $\tilde x(t)$ as a solution to $x'=f(x)+k(t)$ is the difference between the velocity vectors along the path and the underlying vector field, so
\[ k(t)=\frac{d}{dt} \tilde x(t) - f(\tilde x(t)) = f(x_h(t))+h(t)+\frac{v}{T}-f(\tilde x(t)). \]

\noindent We will show that for $\norm{v}$ sufficiently small $\|k-g\|_\infty<r$. Let $L$ be the Lipschitz constant for $f$. For almost every $t\in[0,T]$,
\begin{align*}
\|k(t)-g(t)\|&= \| f(x_h(t))+h(t)+\frac{v}{T}-f(\tilde x(t)) - g(t)  \| \\
&\leq \| f(x_h(t))-f(\tilde x(t))\| + \|h(t) - g(t) \|+ \frac{\norm{v}}{T} \\ 
&\leq L\|x_h(t)-\tilde x(t)\|+  \norm{h(t) - g(t)} + \frac{\norm{v}}{T}\\
&\leq L\norm{v}+ \|h- g \| _\infty+\frac{\norm{v}}{T} \\
&=\norm{v}(L+\frac{1}{T})+ \|h - g \|_\infty
\end{align*}
Let $\epsilon_h=\frac{r-\|h-g\|_\infty}{L+1/T}$. Then $\norm{v}<\epsilon_h$ implies that
$\Infnorm{k-g} < r$. 
Because there exists an $\epsilon_h$-ball about any point in $\cntrlmapp(B_r(g))$ that is also in the image $\cntrlmapp(B_r(g))$, the image of any basis element under $\cntrlmapp$ is open, and the proof is complete. \qed
\end{pf}

With \cref{openG} in hand, we're ready to construct attractor blocks from reachable sets and prove their persistence properties. 

\begin{proposition}\label{stillattrblock}
If $S$ is a nonempty subset of $\R^n$ and if $\reach{r;f}{S}$ is bounded, then for any globally Lipchitz and bounded vector field $\hat{f}:U\rightarrow\R^n$ satisfying $\fnorm{f-\hat{f}}<r$,  the set $\Reach{r;f}{S}$  (defined using $f$) is an attractor block for the flow $\hat\varphi$ generated by $\hat{f}$. 
\end{proposition}

\noindent Before proving \cref{stillattrblock}, we note that in particular it says that boundedness of $\reach{r;f}{S}$ implies that $\Reach{r;f}{S}$ is an attractor block for the flow generated by $f$. Intuitively, one can interpret this result as follows: whenever adding bounded control $g$ results in a bounded reachable set, the underlying vector field $f$ must be counteracting the control at the boundary of the reachable set, ``pulling inwards'' as required for an attractor block.  The strength of this inwards pull reflects the control bound $r$ used to construct the block, and feeds the persistence of the block through autonomous vector field perturbations.

\begin{pf} (\Cref{stillattrblock}) $\Reach{r;f}{S}$ is compact since it is closed and bounded; it is nonempty since it contains $S$. To confirm that $\Reach{r;f}{S}$ is an attractor block for $\hat \varphi$, it will suffice to show that for any $T>0$, 
\begin{equation}
\hat{\varphi}^T(\Reach{r;f}{S}) \overbrace{\subset}^{\text{(I)}} 
\bigcup_{\substack{x_0\in\Reach{r;f}{S} \\ \Infnorm{g}< r}}
\gflowc{T}{x_0}{g} \overbrace{\subset}^{\text{(II)}} \text{int}(\Reach{r;f}{S})
\end{equation}
where the controlled flow in the middle term is generated by the vector field $f$; that is, $\gflowc{T}{x_0}{g}=\gflowf{T}{x_0}{g}{f}$.
Fix $T>0$. For inclusion (I), we first establish that $\hat\varphi^T(x_0)\in \bigcup\limits_{\Infnorm{g}<r} \gflowc{T}{x_0}{g}$ for any point $x_0\in\R^n$. Let the path $\hat x:[0,T]\to\R^n$ be given by $\hat x(t)=\hat\varphi^t(x_0)$; in other words, $\hat x$ solves the initial value problem [$x'=\hat f(x)$, $x(0)=x_0$]. Note that $\hat x$ also solves the initial value problem [$x'=f(x)+g(t)$, $x(0)=x_0$], where $g:[0,T]\rightarrow \R^n$ is given by $g(t)=\hat f(\hat x(t))-f(\hat x(t))$. In the notation of controlled flow for the vector field $f$, $\hat x(T)=\gflowc{T}{x_0}{g}$. Since by hypothesis $\|g\|_\infty\leq \fnorm{f(x)-\hat f(x)}<r$,  it follows that  $\hat\varphi^T(x_0)=\hat x(T)\in \bigcup\limits_{\Infnorm{g}<r} \gflowc{T}{x_0}{g}$. Inclusion (I) follows directly by taking the union over $x_0\in\Reach{r;f}{S}$. 
We argue inclusion (II) by contradiction. Suppose there exists $x^*\in\Reach{r;f}{S}$ and $g^*\in\gspaceI{[0,T]}$, $\Infnorm{g^*}<r$, such that $\gflowc{T}{x^*}{g*}\not\in\text{int}\left( \Reach{r;f}{S}\right)$. Because $\Infnorm{g^*}$ is strictly bounded below $r$, taking $s=\frac{r-\Infnorm{g^*}}{2}$ implies $B_s(g^*)\subset B_r(0)\subset\gspaceI{[0,T]}$. The map $\cntrlmapf{T}{x^*}{f}$ is open (\cref{openG}) so it sends $B_s(g^*)$ to an open neighborhood $\mathcal{O}$ of $\gflowc{T}{x^*}{g^*}$. Because $\gflowc{T}{x^*}{g^*}\not\in\text{int}\left(\Reach{r;f}{S}\right)$, it must be that $\mathcal{O}\cap\compl{\left(\Reach{r;f}{S}\right)}\not=\varnothing$. Therefore, there exists a control $h\in B_s(g)$ such that $\gflowc{T}{x^*}{h}\in\compl{\left(\Reach{r;f}{S}\right)}$, an open set. Continuity of $\gflowc{T}{x}{h}$ with respect to $x$ implies that the inverse image of $\compl{\left(\Reach{r;f}{S}\right)}$ under $\gflowc{T}{\cdot}{h}$ is an open neighborhood $V$ of $x^*$. $V$ must intersect $\reach{r;f}{S}$ nontrivially because $x^*\in\Reach{r;f}{S}$. Let $b\in\reach{r;f}{S}\cap V$, so that $b\in\reach{r;f}{S}$ but $\gflowc{T}{b}{h}\not\in\reach{r;f}{S}$. The contradiction comes from concatenating controls to move from $S$ to $b$ to outside $\reach{r;f}{S}$. In particular, $b\in\reach{r;f}{S}$ implies by \cref{PrS} that there exists $\tau\geq 0$, $a\in S$, and $j\in\gspaceI{[0,\tau]}$, $\Infnorm{j}\leq r$, such that $b=\gflowc{\tau}{a}{j}$. Let $k\in\gspaceI{[0,\tau+T]}$ be given by
\begin{equation*}
k(t)=\begin{cases}
j(t) & \text{ if } 0\leq t <\tau\\
h(t-\tau) & \text{ if } \tau \leq t \leq \tau +T.
\end{cases}
\end{equation*}
Then by construction $\gflowc{\tau+T}{a}{k}=\gflowc{T}{\gflowc{\tau}{a}{j}}{h}=\gflowc{T}{b}{h}\not\in\reach{r;}{S}$. Yet since $\Infnorm{k}\leq r$, \cref{PrS} implies that $\gflowc{\tau+T}{a}{k}\in\reach{r;f}{S}$. This contradiction completes the proof. \qed
\end{pf}

\noindent \Cref{stillattrblock} does not guarantee any relationship between the  set $S$ used to construct the attractor block $\Reach{r;f}{S}$ and the attractor associated with $\Reach{r;f}{S}$. For example, if $f:\R\to\R$ is given by $f(x)=-x(x-1)(x-2)$ and we take a control bound $r=0.5>\sup_{x\in[0,2]} \norm{f(x)}$, then \emph{any} subset $S\subset[0,2]$ will yield the attractor block $\Reach{0.5;f}{S}$ depicted in light blue in \cref{SvA}, with associated attractor $A=\omega(\Reach{r;f}{S})=[0,2]$.
\begin{figure}[h!]
\centering
\includegraphics[scale=0.8]{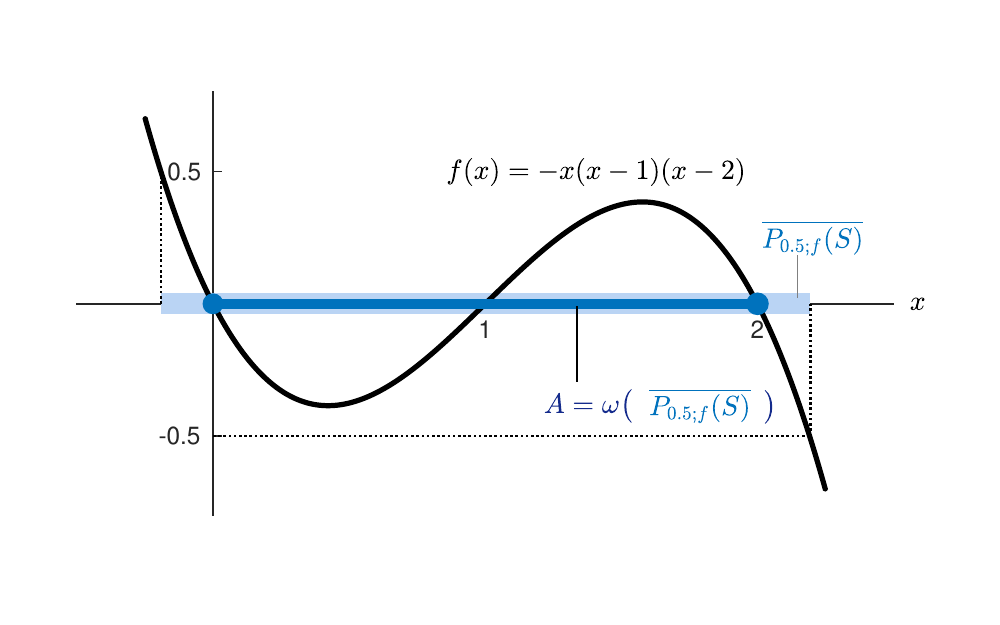}
\vspace{-1cm}
\caption{Relationship between reachable sets and attractors. For $S$ any subset of $[0,2]$, the reachable set under control bound $r=0.5$ is the light blue interval, and the associated attractor is the dark blue interval $[0,2]$.}
\label{SvA}
\end{figure}

On the other hand, if we construct the attractor block $\Reach{r;f}{A}$ from an attractor $A$, we recover $A$ as the associated attractor under one further condition.
\begin{corollary}\label{attrblock}
If $A$ is an attractor and if $\reach{r}{A}$ is contained in a compact subset $K$ of $\domattr{A}$, then $\Reach{r}{A}$ is an attractor block associated with $A$.
\end{corollary}

\begin{pf}\cref{stillattrblock} gives that $\Reach{r}{A}$ is an attractor block. It follows directly from \cref{omega} that the omega limit set preserves inclusion, and so $\omega(A)\subset\omega(\Reach{r}{A})\subset\omega(K)$. By invariance of $A$, $\omega(A)=(A)$. Further, an open cover argument gives that for $K$ a compact subset of $\domattr{A}$, $\omega(K)\subset A$. It follows that $\omega(\Reach{r}{A})=A$, so the attractor block $\Reach{r}{A}$ is associated with the attractor $A$.  \qed
\end{pf}

\subsection{Intensity and Continuation Distance}
\label{dist}
In this subsection we leverage the persistent attractor blocks of \cref{stillattrblock,attrblock} to derive bounds on attractor continuation distance. First, we make precise the notion of continuation.
Suppose that when a vector field $f$ is perturbed to $\hat f$, the original attractor $A$ shifts to another attractor $\hat A$, with qualitatively similar features. In what sense has $A$ continued through the perturbation? We use the following definition, adapted from \cite{mcgehee1988}.
\begin{definition}\label{ic} 
Given an attractor $A$ for vector field $f$ and an attractor $\hat A$ for vector field $\hat f$, we say $A$ \emph{continues immediately} to $\hat A$ if there exists a set $B$ that is an attractor block associated with $A$ under the flow generated by $f$ and is an attractor block associated with $\hat A$ under the flow generated by $\hat f$. 
\end{definition}

\begin{remark}
Because $A$ and its immediate continuation $\hat A$ share a common attractor block, they must also share any topological properties gleaned from the  attractor block.  
\end{remark}

The next theorem gives our central result: an attractor's intensity bounds from below the distance in vector field space over which that attractor continues immediately. 

\begin{theorem}\label{main}
If $A$ is an attractor for vector field $f$ with intensity $\mu(A)$, then for any second vector field $\hat f$ satisfying $\fnorm{f-\hat f} < \mu(A)$, $A$ continues immediately to an attractor $\hat A$ for $\hat f$.
\end{theorem}

\begin{pf}
If $\fnorm{f-\hat f}<\mu(A)$, then there exists a real number $r>0$ such that $\fnorm{f-\hat f}<r<\mu(A)$. It follows from \cref{Idef,rnest2} that $\reach{r}{A}\subset K \subset \domattr{A}$ for some compact set $K\subset\R^n$. This implies, by \cref{attrblock}, that $\Reach{r}{A}$ is an attractor block associated with $A$ for vector field $f$. \Cref{stillattrblock} gives that $\Reach{r}{A}$ is also an attractor block for vector field $\hat f$. Defining $\hat A$ to be the omega limit set of $\Reach{r}{A}$ under $\hat{f}$, we have that $\Reach{r}{A}$ is an attractor block associated with attractors $A$ and $\hat A$ under their respective flows. Hence $A$ continues immediately to $\hat A$, as claimed. \qed
\end{pf}

The following example illustrates continuation of the predator-prey limit cycle from \cref{cycle}.
\begin{example}\label{ppcont}
Consider the predator-prey system \cref{pp1,pp2}, which for suitable parameters features a stable limit cycle $A$ in the first quadrant (see \cref{cyclefig}(c)). Based on reachable set computations with the Euclidean norm, we estimated $0.02 < \mu(A) < 0.03$ in \cref{cycle}. By \cref{main}, $A$ should continue immediately to an attractor $\hat{A}$ for any second vector field that differs from the original by up to $0.02$ in Euclidean sup-norm. \Cref{cyclepert} illustrates the immediate continuation $\hat{A}$ for three such vector field perturbations: (a) adding 0.02 to $dx/dt$, (b) subtracting 0.02 from $dy/dt$, and (c) changing the prey carrying capacity parameter from $K=4$ to $K=0.39801$ (this parameter change perturbs the vector field by less than 0.02 in a neighborhood of the attractor block $\Reach{0.02}{A}$). The original cycle $A$ is shown with a solid line, the perturbed attractor $\hat{A}$ is shown with a dashed line, and their common attractor block $\Reach{0.02}{A}$ is shaded in green. \hfill //
\end{example}

\begin{figure}[h]
\centering
\includegraphics[width=\textwidth]{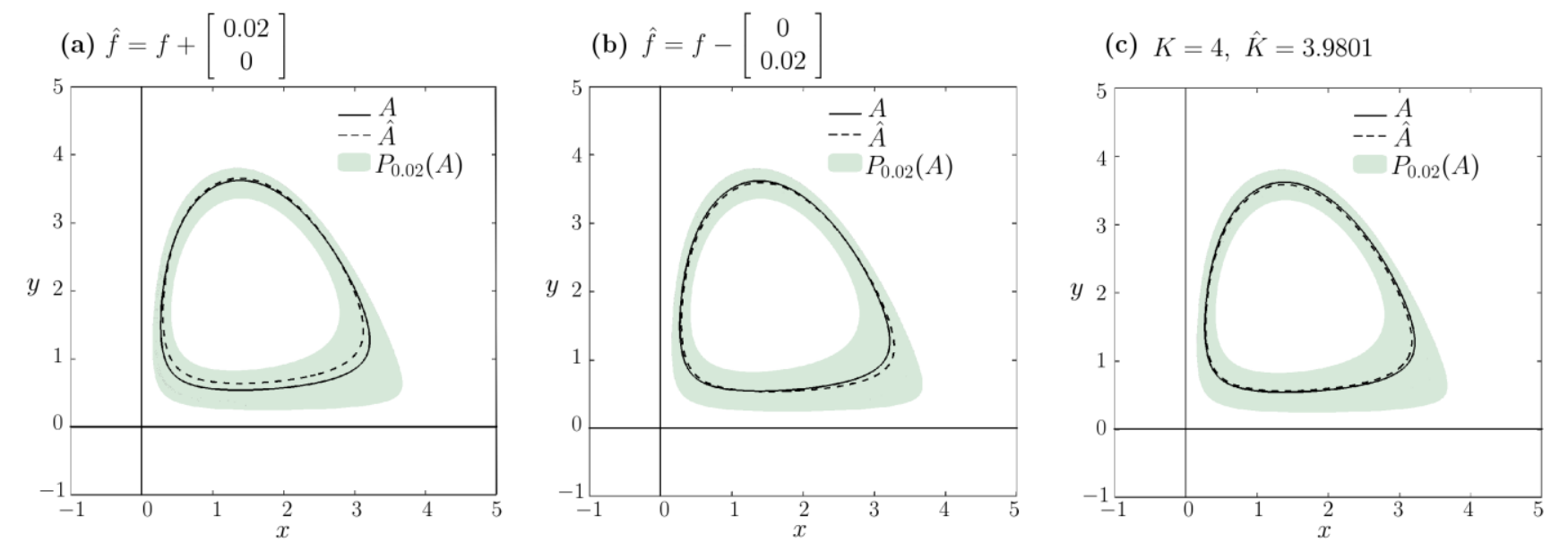}
\vspace{-1cm}
\caption{Continuation of a predator-prey limit cycle through three time-independent pertubations to the vector field $f$. The perturbed vector field $\hat f$ is given above each panel. Solid and dashed lines give the original and perturbed limit cycle, respectively. The reachable set from the original limit cycle with control bound $r=0.02$ is shaded in green.\\}
\label{cyclepert}
\end{figure}

\Cref{ppcont} highlights the strength of \cref{main} in guaranteeing attractor continuation across a variety of perturbation forms, rather than perturbation of a specific parameter. 

The closing example in this subsection demonstrates that attractors may continue for vector field perturbations that \emph{exceed} intensity of attraction. This does not contradict \cref{main}, but indicates that intensity does not give an upper bound on attractor continuation distance.
\begin{example}\label{notub}
Let $x'=f(x)=x(1-x)$.  Then $A=1$ is an attractor for $f$ with intensity $\mu(A)=0.25$. The set $B=\left[\frac{1}{2},\frac{3}{2} \right]$ is an attractor block associated with $A$. Consider a second system $x'=\hat f(x)=cx(1-x)$ with $c>0$. $B$ remains an attractor block associated with $\{1\}$ for any such system. Yet the distance $\fnorm{f-\hat f}$ may exceed $0.25$ by taking $c$ sufficiently large. Thus the attractor $A=1$ continues immediately to $\hat A=1$ despite $\fnorm{f-\hat f}$ exceeding $\mu(A)$.  \hfill //
\end{example}

\subsection{Upper Semicontinuity of Attractors}
\label{upper}
Without too much extra effort, the present framework yields a new proof that attractors are upper semicontinuous---roughly, they do not expand suddenly in response to changes in a vector field. Formally,
\begin{definition}\label{usc}
An attractor $A$ corresponding to the vector field $f$ is \emph{upper semicontinuous} if for any neighborhood $V$ of $A$ there exists a $\delta>0$ such that any for any second vector field $\hat{f}$ with $\fnorm{f-\hat{f}}<\delta$, $A$ continues to an attractor $\hat{A}\subset V$. 
\end{definition}

In this subsection we clear some debts by proving \cref{PrAtoA,mug0}, then combine \cref{PrAtoA} with results on attractor continuation to prove upper semicontinuity of attractors in our setting (\cref{uppersc}).

\begin{pf}
(\cref{PrAtoA})
Property (1) is a direct instance of \cref{rnest2}. Towards property (2), fix a neighborhood $V$ of $A$. \Cref{ab} gives existence of an attractor block $B\subset \text{int}(V)$ associated with $A$. Let $\delta_1=\text{dist}(B,\overline{V^C})$ and $\delta_2=\text{dist}(\varphi^1(B),\overline{B^C})$. Let $\delta=\text{min}\{\delta_1,\delta_2\}$  and note $\delta>0$. By \cref{contg} there exists an $\hat r>0$ such that for all $x\in\R^n$ and all $t\in[0,1]$, $\Infnorm{g}<\hat r$ ensures $\norm{\gflowc{t}{x}{g}-\varphi^t(x)}<\delta/2$. Let $r=\hat r /2$. Then $\Infnorm{g}\leq r$ implies that for $t\in[0,1]$,
\begin{equation}\label{s}
\gflowc{t}{B}{g}\subset\mathcal{N}_{\delta/2}(\varphi^t(B))\subset \mathcal{N}_{\delta/2}(B)\subset V
\end{equation}
where the second inclusion follows from forward invariance of $B$ and the third from the construction of $\delta$. Additionally, at time $1$ we have
\begin{equation}\label{tau}
\gflowc{1}{B}{g}\subset\mathcal{N}_{\delta/2}(\varphi^1(B))\subset B
\end{equation}
with the final inclusion again following from construction of $\delta$.
A simple inductive argument based on inclusions \cref{s} and \cref{tau} implies that $\gflowc{t}{B}{g}\subset V$ for all $t\geq 0$. Since $A\subset B$,  this implies $\gflowc{t}{A}{g}\subset V$ for all $t\geq 0$ and all $\Infnorm{g}\leq r$. By \cref{PrS} we have $\reach{r}{A}\subset V$, establishing property (2).
The inclusion $A\subset\bigcap_{r>0}\reach{r}{A}$ of property (3) is immediate. The reverse inclusion follows from property (2): any point not in $A$ can be excluded from some neighborhood of $A$ and hence from $\reach{r}{A}$ for some $r>0$. Hence any point in $\reach{r}{A}$ for all $r>0$ must be in $A$. This completes the proof. \qed
\end{pf}

\begin{pf}
 (\cref{mug0}) We wish to show that $\mu(A)>0$ for any attractor $A$. \Cref{ab} gives an attractor block $B$ associated with $A$. One can readily verify that $B\subset\mathcal{D}(A)$.
By \cref{PrAtoA}, there exists an $r>0$ such that $\reach{r}{A}\subset B \subset \mathcal{D}(A)$. \Cref{Idef} then implies that $\mu(A)>0$.  \qed
\end{pf}

\begin{theorem}
\label{uppersc}
Attractors in the present setting are upper semicontinuous in the sense of \cref{usc}.
\end{theorem}

\begin{pf}
Fix a neighborhood V of an attractor $A$ for $x'=f(x)$. We will show that there exists a $\delta$ such that $\fnorm{f-\hat{f}}<\delta$ implies that $A$ continues immediately to an attractor $\hat{A}\subset V$ for $\hat{f}$. Using \cref{compactnbhd}, let $K$ be a compact set with $A\subset \text{int}(K)\subset K \subset V\cap \domattr{A}$. Part 2 of \cref{PrAtoA} gives an $r>0$ such that $\reach{r}{A}\subset\text{int}(K)$. Then $\Reach{r}{A}\subset K \subset \domattr{A}$. Let $\delta=r$. By \cref{stillattrblock}, $\fnorm{f - \hat f} < r$ implies that $\Reach{r}{A}$ is an attractor block for $\hat f$. Let $\hat A$ be the omega limit set of $\Reach{r}{A}$ under $\hat f$. Then $\hat A$ is an attractor under $\hat{f}$ (\cref{obattr}), is an immediate continuation of $A$ (\cref{ic}), and is contained in V, as desired. \qed
\end{pf}

\section{Conclusions}

We have defined a quantity, intensity of attraction, which gives a lower bound on the magnitude of vector field perturbations through which an attractor continues to structure the long-term behavior of trajectories. Intensity is an inherently metric property of an attractor, and its value depends on the metric used on $\R^n$. Provided that time-dependent perturbations to a vector field are essentially bounded below the intensity of an attractor, trajectories that start at the attractor remain within its domain of attraction. And when one autonomous vector field is perturbed to a second one, an attractor for the first continues immediately to the second as long as the time-independent perturbation does not meet or exceed its intensity.  

An attractor and its immediate continuation share the topological properties encoded in their common attractor block---and these properties are generally more coarse than the details treated in earlier works on structural stability \cite{andronov1937rough,peixoto1962structural,smale1967differentiable}. By accepting this coarser lens on system structure we get attractor persistence not just for ``sufficiently small'' vector field perturbations, as in \cite{conleyeaston1971, conley1978}, but specifically for perturbations smaller than the attractor's intensity. 

We envision several extensions of the present work in both applied and theoretical directions. Efficient numerical algorithms will be an important bridge between intensity theory and the real-world modeling problems that motivated its development. To compute the intensity of an attractor we need information about both the domain of attraction and a nested collection of sets reachable from the attractor. Delineating a domain of attraction in a nonlinear system can be a non-trivial problem in its own right, and a variety of computational approaches have been developed (see \cite{matallana2010estimation} and references therein). To compute reachable sets quickly, one could pursue  numerical strategies described in \cite{sethian1999level}. To achieve rigorous bounds on reachable sets, one could use methods based on interval arithmetic (\cite{deWeerdt2011}).

Applying intensity theory to a real-world system will require selecting a metric that best reflects the perturbations of interest. Here we have considered homogenous, translation-invariant metrics, carrying the same information as a norm. However, the relevant scale of perturbations might vary across state space, requiring a metric in its full generality. Furthermore, in certain settings one may want to restrict perturbations to a single coordinate while leaving others unperturbed. This could be accomplished with an extended metric. Whether intensity theory goes through for general or extended metrics remains to be determined. 

Another theoretical question with practical implications concerns the relationship between intensity of attraction for discrete \cite{mcgehee1988} and continuous-time dynamics. We conjecture that $\mu(A)=\lim\limits_{t\to 0} \mu_t (A) / t$, where $\mu$ gives intensity of an attractor under a flow $\varphi$ (\cref{Idef}) and $\mu_t$ gives intensity of the attractor under the map $\varphi^t$ (defined in \cite{mcgehee1988}, section 5). Proving this connection would provide a theoretical justification for discretizing time in numerical computations of intensity $\mu$.

Lastly, we aim to extend the theory of intensity to measure persistence properties of repellers and other isolated invariant sets.  If we can achieve bounds on the persistence of local dynamic structures, it might be possible to describe how the global structure changes, in stages, under perturbations of increasing magnitude. We conjecture that this series of changes will reveal a Morse decomposition of the dynamical system.

\section*{Acknowledgements}
\noindent This work was funded by the National Science Foundation [grant numbers 00039202, DMS-1645643, DMS-0940366, DMS-0940363]. We thank members of the McGehee research group, the Math Climate Research Network, and the Cornell University Mathematics Department who helped push this work forward. In particular, thanks go to Alice Nadeau, Shannon Negaard Paper, Bill Satzer, Mary Lou Zeeman, Kelly Patwell, Alex Vladimirksy, John Hubbard, and Steven Strogatz for helpful conversations and support. Cornell University supported K. Meyer during the writing phase of this project, and Steven Strogatz and Mary Lou Zeeman provided feedback on an earlier draft of the introduction of the manuscript.

\bibliographystyle{siam}
\bibliography{MPAbib}

\end{document}